# A remarkable periodic solution of the three-body problem in the case of equal masses

By Alain Chenciner and Richard Montgomery

*Dedicated to Don Saari for his (censored) birthday*

### Abstract

Using a variational method, we exhibit a surprisingly simple periodic orbit for the newtonian problem of three equal masses in the plane. The orbit has zero angular momentum and a very rich symmetry pattern. Its most surprising feature is that the three bodies chase each other around a fixed eight-shaped curve. Setting aside collinear motions, the only other known motion along a fixed curve in the inertial plane is the "Lagrange relative equilibrium" in which the three bodies form a rigid equilateral triangle which rotates at constant angular velocity within its circumscribing circle. Our orbit visits in turns every "Euler configuration" in which one of the bodies sits at the midpoint of the segment defined by the other two (Figure 1). Numerical computations by Carles Simó, to be published elsewhere, indicate that the orbit is "stable" (i.e. completely elliptic with torsion). Moreover, they show that the moment of inertia $I(t)$ with respect to the center of mass and the potential $U(t)$ as functions of time are almost constant.

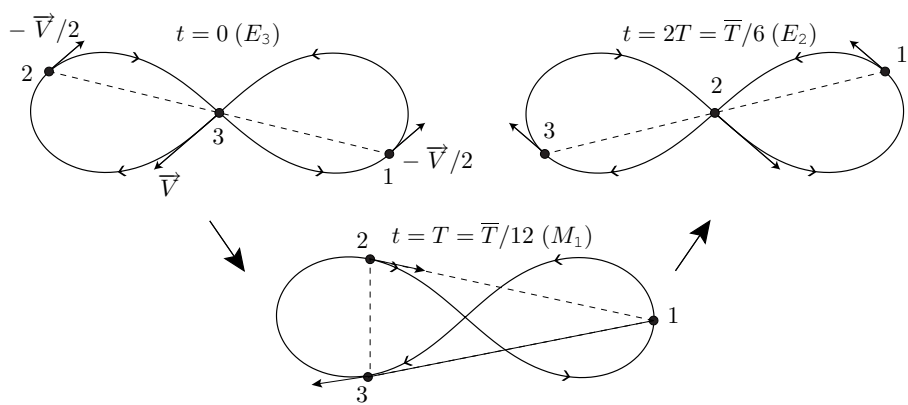

Figure 1 (Initial conditions computed by Carles Simó)
$x_1=-x_2=0.97000436-0.24308753i, x_3=0$; $\vec{V}=\dot{x}_3=-2\dot{x}_1=-2\dot{x}_2=-0.93240737-0.86473146i$
$\overline{T}=12T=6.32591398$, $I(0)=2$, $m_1=m_2=m_3=1$

## 1. The setting

We consider three bodies of unit mass in the Euclidean plane $\mathbb{R}^2 = \mathbb{C}$, with Newtonian attraction. After using Galilean invariance to fix the center of mass, the phase space becomes the tangent bundle to the configuration space $\hat{\mathcal{X}} = \mathcal{X} \setminus \{\text{collisions}\}$,

$$\mathcal{X} = \left\{ x = (x_1, x_2, x_3) \in (\mathbb{C})^3, \sum_{i=1}^{3} x_i = 0 \right\}.$$

We endow $\mathcal{X}$ with the *Hermitian mass scalar product*. In the case of equal unit masses this is

$$\langle x, y \rangle = \sum_{j=1}^{3} \overline{x}_j \cdot y_j = x \cdot y + i\omega(x, y).$$

Its real and imaginary part are respectively the *mass scalar product* and the *mass symplectic structure*. The isometry group $O(2)$ of $\mathbb{R}^2 = \mathbb{C}$ acts diagonally on $\mathcal{X}$: the symmetry $S$ with respect to the first coordinate axis and the rotation $R_\theta$ of angle $\theta$ act respectively as

$$\begin{cases} S \cdot (x_1, x_2, x_3) = (\bar{x}_1, \bar{x}_2, \bar{x}_3), \\ R_\theta \cdot (x_1, x_2, x_3) = (e^{i\theta} x_1, e^{i\theta} x_2, e^{i\theta} x_3). \end{cases}$$

The phase space will be identified with the Cartesian product $\hat{\mathcal{X}} \times \mathcal{X}$ with elements written $(x, y)$. We define the following $O(2)$-invariant functions on phase space:

$$I = x \cdot x, \ J = x \cdot y, \ K = y \cdot y, \ U = U(x), \ H = \frac{1}{2}K - U, \ L = \frac{1}{2}K + U.$$

These are the *moment of inertia with respect to the center of mass*, half its derivative with respect to time, twice the *kinetic energy* in a frame attached to the center of mass, the *potential function*, the *total energy*, and the *Lagrangian*. The *size* $r = I^{\frac{1}{2}}$ of the configurations is the norm on $\mathcal{X}$. The "force function" $U$, the negative of the potential energy, is defined by

$$U = \frac{1}{r_{12}} + \frac{1}{r_{13}} + \frac{1}{r_{23}}, \quad \text{where} \quad r_{ij} = |x_i - x_j|.$$

(We choose units so that the gravitational constant is 1.)

The *central configurations* play a basic role in this paper. They are the only configurations which admit *homothetic motions*, that is motions for which the configuration collapses homothetically to its center of mass. They admit more generally *homographic motions* where each body has a similar Keplerian motion and in particular *relative equilibrium motions* where the bodies rotate rigidly and uniformly around their center of mass. These configurations are the critical points of the *scaled potential function* $\tilde{U} = \sqrt{I}U$. Upon normalization these are the critical points of the restriction $U|_{I=1}$ of the potential function to the sphere $I = 1$. In the case of three bodies, the central configurations are completely known thanks to the works of Euler and Lagrange. There are three collinear configurations $E_1, E_2, E_3$, distinguished by which mass sits in between the two others (it sits at the midpoint in the case of equal masses), and two equilateral ones $L^+, L^-$ distinguished only by their orientation.

In naming the central configurations we have already formed the quotient of the configuration space by the rotation group SO(2). It is well-known ([6], [7]) that after reduction by direct isometries (translations and rotations) the three-body problem in the plane has a configuration space which is homeomorphic to $\mathbb{R}^3$. This reduced configuration space – the space of oriented triangles in the plane up to translation and rotation – is endowed with a metric induced from the mass metric on configuration space which makes it a cone over a round 2-sphere of radius $\frac{1}{2}$. Points of this sphere, henceforth called the *shape sphere*, represent oriented triangles whose moment of inertia $I$ is 1. This sphere is to be thought of as the space of oriented similarity classes of triangles [7]. The shape sphere is depicted in Figure 2 in the case of equal masses; it will be studied from a Riemannian point of view in paragraph 5. Its main features are the North and South poles, corresponding to the the Lagrange points $L^+$ and $L^-$, the equator corresponding to collinear triangles and the three meridians corresponding to the three types of isosceles triangles. Each meridian $M_i$ intersects the equator at an Euler point $E_i$ and an antipodal collision point $C_i$. The figure also shows the equipotential curve (level curve of $U|_{I=1}$) containing the Euler points and a level surface of $U$ corresponding to a higher value. The projection of the "eight" orbit onto the shape sphere closely resembles the equipotential curve passing through the three Euler points, and shares its symmetries. We have stressed the half-line $E_3$ and the meridian plane $M_1$ because they explicitly appear in the proof. Note that we used a false perspective for the poles to increase readability.

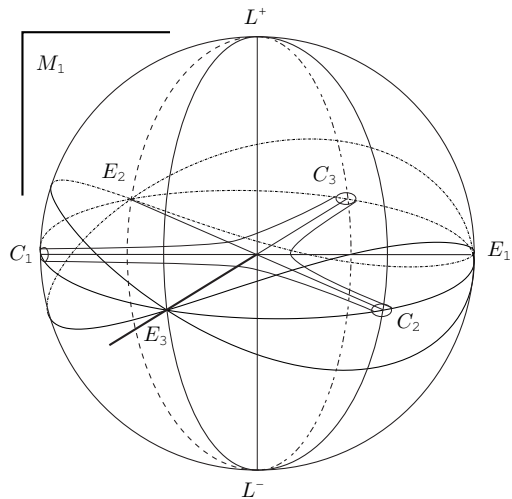

Figure 2

## 2. The orbit

Let $\overline{T}$ be any positive real number. We define actions of the Klein group $\mathbb{Z}/2\mathbb{Z} \times \mathbb{Z}/2\mathbb{Z}$ on $\mathbb{R}/\overline{T}\mathbb{Z}$ and on $\mathbb{R}^2$ as follows: if $\sigma$ and $\tau$ are generators,

$$\sigma \cdot t = t + \frac{\overline{T}}{2}, \quad \tau \cdot t = -t + \frac{\overline{T}}{2}, \quad \sigma \cdot (x,y) = (-x, y), \quad \tau \cdot (x,y) = (x, -y).$$

THEOREM. *There exists an "eight"-shaped planar loop $q : (\mathbb{R}/\overline{T}\mathbb{Z}), 0 \to \mathbb{R}^2, 0$ with the following properties*:

(i) *for each $t$,*
$$q(t) + q(t + \overline{T}/3) + q(t + 2\overline{T}/3) = 0;$$

(ii) *$q$ is equivariant with respect to the actions of $\mathbb{Z}/2\mathbb{Z} \times \mathbb{Z}/2\mathbb{Z}$ on $\mathbb{R}/\overline{T}\mathbb{Z}$ and $\mathbb{R}^2$ defined above*:
$$q(\sigma \cdot t) = \sigma \cdot q(t) \quad \text{and} \quad q(\tau \cdot t) = \tau \cdot q(t);$$

(iii) *the loop $x : \mathbb{R}/\overline{T}\mathbb{Z} \to \hat{\mathcal{X}}$ defined by*
$$x(t) = \big(q(t + 2\overline{T}/3), q(t + \overline{T}/3), q(t)\big)$$

*is a zero angular momentum $\overline{T}$-periodic solution of the planar three-body problem with equal masses.*

The rest of the paper is devoted to a proof of this theorem.

## 3. The structure of the proof

The main ingredients of the proof are the following:

(i) *The "direct" method.* One twelfth of the orbit is obtained by minimizing the action (we note $T = \overline{T}/12$)

$$\mathcal{A} = \int_0^T (\frac{1}{2}K + U)\,dt$$

over the subspace $\Lambda$ of $H^1([0,T], \mathcal{X})$ consisting of paths which start at an Euler configuration, say $E_3$ (3 in the middle) with arbitrary size and end at an isosceles configuration, say $M_1$ ($r_{12} = r_{13}$), again of arbitrary size. Existence of a minimizing path is standard. The main point is to show that such a path has no collision.

(ii) *Reduction.* The kinetic energy can be expressed as the sum of two nonnegative terms: $K = K_{\text{red}} + K_{\text{rot}}$. (This is Saari's decomposition of the velocity. See paragraph 5 below.) $K_{\text{red}}$ corresponds to a Riemannian metric on the quotient space $\mathcal{X}/\text{SO}(2)$ induced from the metric $K$ on $\mathcal{X}$. It is the deformation part of the kinetic energy, including the homothetic part. (See Lemma 5 below for its explicit expression.) $K_{\text{rot}}$ is the rotational part of kinetic energy. Because we consider the planar problem, $K_{\text{rot}} = |\mathcal{C}|^2/I$ with $\mathcal{C}$ the angular momentum vector. (This equality is replaced by the inequality $K_{\text{rot}} \geq |\mathcal{C}|^2/I$ in the spatial three-body problem. See Sundman's inequality $IK - J^2 \geq 0$ in [1].) The boundary conditions defining $\Lambda$ are invariant under rotation, which means that if $x(t)$ is such a path, so is $g(t)x(t)$ for $g(t)$ any sufficiently regular curve in SO(2). Changing the choice of $g(t)$ changes $K_{\text{rot}}$ while leaving $K_{\text{red}}$ fixed. By appropriate choice of $g(t)$ we can insure that $K_{\text{rot}} = 0$ and hence $\mathcal{C} = 0$. This shows that any minimizer for the problem of (i) has zero angular momentum. Moreover, the original problem has been reduced to the problem of minimizing the reduced action

$$\mathcal{A}_{\text{red}} = \int_0^T \left(\frac{1}{2}K_{\text{red}} + U\right) dt$$

over paths lying in the reduced configuration space $\mathcal{X}/\text{SO}(2)$, and satisfying the same boundary conditions. The advantage of this reduction is that the reduced configuration space is topologically $\mathbb{R}^3$, hence paths are easy to visualize (see Figures 2, 3, and 4).

(iii) *Comparison with Kepler to exclude collisions.* Instead of computing local variations of the action, we compute the infimum of the action over all paths in $H^1([0,T], \mathcal{X})$ with collision. Denote this infimum by $A_2$, the "2" standing for two bodies because we explicitly compute $A_2$ via a two body problem.

We compare $A_2$ to the action $a$ of an explicit carefully chosen collision-free "test path" in $\Lambda$. We show that $A_2 > a$, yielding the result that the minimizer must be collision-free.

(iv) *Symmetries and area rule.* The equality of masses gives a supply of symmetries corresponding to interchanging the masses. Using these, we construct eleven other congruent copies of the minimizer from (i). The first variation of action formula shows that these copies fit together smoothly with the original to form a single orbit, periodic in the reduced configuration space (i.e. modulo rotations), of period $12T$. To reconstruct the motion in the inertial plane, i.e. in $\mathcal{X}$, we use the symmetries combined with the "area formula" (see [9]). These tools yield that the motion in $\mathcal{X}$ is periodic, that all three masses indeed move along one and the same curve in the inertial plane, with the Klein group symmetry described in the theorem.

(v) *Proving the curve has the shape of a "figure eight".* We prove that the individual angular momentum of any one of the three bodies is zero only as this body passes through the origin. This implies that each one of the two lobes of the curve is starshaped.

*Other proofs.* We have presented the shortest proof we know. Another proof proceeds by local perturbations so as to destroy collisions while decreasing the action, getting rid of the binary and triple collisions separately. Special care must be taken with the direction of perturbation used when destroying binary collisions. Estimates are significantly longer than those given here, but do not require numerical integration. On the other hand, they do not yield the estimates on the orbit given in Appendix 1.

## 4. The exclusion of collisions

Let us recall that $\Lambda$ is the subspace of $H^1([0,T], \mathcal{X})$ consisting of paths which start at the Euler configuration $E_3$ with arbitrary size and end at an isosceles configuration of type $M_1$ with arbitrary size. This paragraph and the next one are devoted to proving

PROPOSITION. *A path in $\Lambda$ which minimizes the action*

$$\mathcal{A} = \int_0^T (\frac{1}{2}K + U)\,dt$$

*has no collisions.*

*Proof of the proposition.* Surprisingly, we are able to treat double and triple collisions simultaneously. The first key point is the (trivial) remark that

the action

$$\mathcal{A}(m_1, m_2, m_3; x) = \int_0^T \left( \frac{1}{2} \sum_{i=1}^{3} m_i |\dot{x}_i(t)|^2 + \sum_{1 \leq i < j \leq 3} \frac{m_i m_j}{|x_j(t) - x_i(t)|} \right) dt$$

along a given path $t \mapsto x(t) = \bigl(x_1(t), x_2(t), x_3(t)\bigr)$ is an increasing function of any of the masses. Setting, for example, $m_1$ to zero, yields

$$\mathcal{A}(x) = \mathcal{A}(1, 1, 1; x) > \mathcal{A}(0, 1, 1; x).$$

The last term is the action of a 2-body problem with equal masses. We use this remark in the following lemma:

LEMMA 1. *If $x \in H^1([0, T], \mathcal{X})$ has a collision, either double or triple, then its action is greater than $A_2$, where $A_2$ is the action of the two-body collinear solution in which two unit masses start at collision and end with zero velocity at a time $T$ later, the center of mass being fixed throughout. (This might be called half a collision-ejection elliptic orbit.) See below for an explicit formula for $A_2$.*

*Remark.* The lemma asserts that the infimum of the action $\mathcal{A}(x)$ over all collision paths $x \in H^1([0, T], \mathcal{X})$ is greater than or equal to $A_2$. Indeed this infimum equals $A_2$. Imagine a sequence $x_n$ of paths in which $m_2$ and $m_3$ perform half the Keplerian collision-ejection orbit described in the lemma, while $m_1$ remains fixed, and at distance $n$ from the 2-3 center of mass. We have $\mathcal{A}(x_n) \to A_2$ as $n \to \infty$.

*Proof of Lemma* 1. Let us suppose that masses 2 and 3 (and possibly 1) collide at the instant $T_1$. As just explained, we lower the path's action by setting $m_1$ equal to zero. We can forget about the position of mass 1 for this new action, and are left with the action for the Keplerian two body problem, investigated by Gordon in [4]. According to Gordon each part of the curve $x$, before $T_1$ and after $T_1$, has action greater than or equal to the corresponding collinear motion of masses $m_2$ and $m_3$ in which they collide at $T_1$ and are at rest at their other endpoints, $t = 0$ or $t = T$. Indeed, by doubling each part by concatenating it with itself but reversed, we get two closed paths, each going from collision to collision, one in time $2T_1$ and the other in time $2(T - T_1)$. The absolute minimum for the collision-to-collision problem was shown by Gordon to be the collision-ejection solution. Its action is proportional to $T^{1/3}$, a convex function of the period $T$. Following Gordon, this convexity implies that the action is lowered further if we replace our previous two motions by the single collision-ejection solution beginning and ending at collision with no collisions in between. Half of this path realizes the infimum to collision in time $T$ for the Kepler problem. This ends the proof.

The next lemma introduces our collisionless test path and reduces the proof of the proposition to an estimation of the length $\ell_0$ of its projection on the shape sphere.

*Equipotential test paths.* Fix $I = I_0$ and $U = U_0$ where $U_0$ is the value of the potential function at any one of the Euler configurations on the sphere $I = I_0$. Viewed in the reduced configuration space, this defines a curve on the two-sphere of radius $\sqrt{I_0}$ (see Figure 2). Take the one-twelfth of this curve lying above the equator, connecting $E_3$ to $M_1$. Traverse this curve at constant speed, the speed chosen so as to finish at the desired time $T$. This gives us a family of reduced test paths in $\mathcal{X}/\mathrm{SO}(2)$ depending on $I_0$. The corresponding paths in $\mathcal{X}$ are those which have zero-angular momentum and project to these. The lengths of these paths are $\ell_0\sqrt{I_0}$ where $\ell_0$ is the length of the path when $I_0 = 1$, henceforth called the "Euler equipotential length."

LEMMA 2. *The minimum $a$ of the actions of equipotential test paths is smaller than $A_2$, the infimum of the actions of collision paths in time $T$, if and only if the Euler equipotential length $\ell_0$ satisfies*

$$\ell_0 < \frac{\pi}{5}.$$

*Proof of Lemma* 2. The action $A_2$ of Lemma 1 is half that of the collision-ejection path of period $2T$. It is computed (for example in [2] ; take $\kappa = -\frac{1}{2}$ there) to be:

$$A_2 = \frac{1}{2} \times 3 \times (2\pi^2)^{\frac{1}{3}}(\frac{1}{2}\tilde{U}_2)^{\frac{2}{3}}(2T)^{\frac{1}{3}}; \quad \tilde{U}_2 = \frac{1}{\sqrt{2}}.$$

The constant $\tilde{U}_2$, which is called $U_0$ in [2], is the (constant) value of the scaled potential function $\tilde{U} = \sqrt{I}U$ for the 2-body problem with both masses equal to 1.

We now evaluate the minimum $a$ of the actions of the equipotential test paths. (This computation is the same as the one in [2].) In computing the action $A(I_0)$ of the test path at radius $\sqrt{I_0}$ note that both integrands are constant. The action is then

$$A(I_0) = (\frac{1}{2}K_0 + U_0)T, \quad \text{where} \quad K_0 = \left(\frac{\ell_0\sqrt{I_0}}{T}\right)^2, \quad U_0 = \frac{\tilde{U}_E}{\sqrt{I_0}}.$$

The constant $\tilde{U}_E = \frac{5}{\sqrt{2}}$ is the value of $\tilde{U} = \sqrt{I}U$ at the Euler configurations. As in [2], we are left with the problem

$$\text{minimize: } A(I_0) = \frac{1}{2}\left(\frac{\ell_0\sqrt{I_0}}{T}\right)^2 T + \frac{5}{\sqrt{2}}\frac{1}{\sqrt{I_0}}T.$$

This function has a unique minimum with respect to $I_0$ at

$$I_0 = \left(\frac{5}{\sqrt{2}\ell_0^2}\right)^{\frac{2}{3}} T^{\frac{4}{3}}.$$

The corresponding action is

$$a = \frac{3}{2}\left(\frac{5}{\sqrt{2}}\right)^{\frac{2}{3}} \ell_0^{\frac{2}{3}} T^{\frac{1}{3}}.$$

Finally, $a < A_2$ if and only if

$$\frac{3}{2}\left(\frac{5}{\sqrt{2}}\right)^{\frac{2}{3}} \ell_0^{\frac{2}{3}} T^{\frac{1}{3}} < \frac{1}{2} \times 3 \times (2\pi^2)^{\frac{1}{3}}(\frac{1}{2}\tilde{U}_2)^{\frac{2}{3}}(2T)^{\frac{1}{3}},$$

which amounts to

$$\ell_0 < \frac{\sqrt{2}}{5}\tilde{U}_2\,\pi = \frac{\pi}{5}.$$

## 5. Length computations

LEMMA 3. *The Euler equipotential length satisfies $\ell_0 < \pi/5$.*

In order to get this estimate we will need explicit coordinates on the quotient space $\mathcal{X}/\mathrm{SO}(2)$, and expressions for the metric and the potential function in these coordinates. The estimation of $\ell_0$ is then obtained numerically with great accuracy.

(i) *The quotient map.* One realizes the quotient by composing the "Jacobi map" $\mathcal{J}$ with the "Hopf map" $\mathcal{K}$. The configuration space $\mathcal{X}$ is a two-dimensional complex Hermitian vector space. As such, it is isometric to $\mathbb{C}^2$. Jacobi coordinates

$$\mathcal{J}: \mathcal{X} \to \mathbb{C}^2$$

defined by

$$(z_1, z_2) = \mathcal{J}(x_1, x_2, x_3) = \left(\frac{1}{\sqrt{2}}(x_3 - x_2), \sqrt{\frac{2}{3}}\left(x_1 - \frac{1}{2}(x_2 + x_3)\right)\right)$$

realize this isomorphism. (Recall that if $(x_1, x_2, x_3)$ is a point in $\mathcal{X}$ then $x_1 + x_2 + x_3 = 0$.) Being an isometry, we have $I = |z_1|^2 + |z_2|^2$ in Jacobi coordinates. The action of $\mathrm{SO}(2)$ corresponds to the diagonal action of the complex unit scalars on Jacobi coordinates $(z_1, z_2)$. That is to say, $\mathcal{X}/\mathrm{SO}(2) = \mathbb{C}^2/S^1$. The quotient by rotations is realized by making a vector out of the invariant polynomials for this action:

$$\mathcal{K}(z_1, z_2) = (u_1, u_2 + iu_3) = \left(|z_1|^2 - |z_2|^2, 2\bar{z}_1 z_2\right).$$

This is the Hopf (also called Kustaanheimo-Stiefel) map

$$\mathcal{K} : \mathbb{C}^2 \to \mathbb{R} \times \mathbb{C} = \mathbb{R}^3.$$

*Remark.* A compact definition of $\mathcal{K}$ is obtained by identifying $\mathbb{C}^2$ (respectively $\mathbb{R}^3$) with the quaternions $\mathbb{H}$ (respectively with the purely imaginary quaternions $\tilde{\mathbb{H}}$) as follows:

$$\mathbb{C}^2 \ni (z_1, z_2) \mapsto q = z_1 + z_2 j \in \mathbb{H}$$
$$\mathbb{R}^3 \ni (u_1, u_2, u_3) \mapsto u_1 i + (u_2 + i u_3) k = u_1 - u_3 j + u_2 k \in \tilde{\mathbb{H}}.$$

Then $\mathcal{K} : \mathbb{H} \to \tilde{\mathbb{H}}$ is defined by $\mathcal{K}(q) = \bar{q} i q$.

Here are some properties of these mappings:

$$|\mathcal{K}(z_1, z_2)|^2 = u_1^2 + u_2^2 + u_3^2 = (|z_1|^2 + |z_2|^2)^2 = I^2.$$

The 3-sphere $I = 1$ is sent by $\mathcal{K} \circ \mathcal{J}$ to the unit 2-sphere of $\mathbb{R}^3$ (the *shape sphere*) according to the *Hopf fibration*. Indeed, restricted to $I = 1$, our formula for $\mathcal{K}$ is the standard formula for the Hopf fibration. The location on this sphere of the collision points, the Euler points and the Lagrange points are:

$$C_1 = (-1, 0, 0), \quad C_2 = (\frac{1}{2}, \frac{\sqrt{3}}{2}, 0), \quad C_3 = (\frac{1}{2}, -\frac{\sqrt{3}}{2}, 0),$$
$$E_1 = (1, 0, 0), \quad E_2 = (-\frac{1}{2}, -\frac{\sqrt{3}}{2}, 0), \quad E_3 = (-\frac{1}{2}, \frac{\sqrt{3}}{2}, 0),$$
$$L^+ = (0, 0, 1), \quad L^- = (0, 0, -1).$$

Using the above formulas and the expressions

$$r_{23} = \sqrt{2}|z_1|, \quad r_{31} = |\sqrt{3/2} z_2 + (1/\sqrt{2}) z_1|, \quad r_{12} = |\sqrt{3/2} z_2 - (1/\sqrt{2}) z_1|,$$

the proof of the following lemma is immediate:

LEMMA 4 (Hsiang [5]). *To $u = (u_1, u_2, u_3)$ in the shape sphere, corresponds a triangle with sides $r_{23} = \sqrt{1 - C_1 \cdot u}$, $r_{31} = \sqrt{1 - C_2 \cdot u}$, $r_{12} = \sqrt{1 - C_3 \cdot u}$, where the scalar product is the standard Euclidean one in $\mathbb{R}^3$.*

(ii) *The orbit metric.* We derive a formula for the reduced metric $K_{\text{red}}$ by computing the distance $d(x, y)$ between the SO(2)-orbits of $x$ and $y$ in $\mathcal{X}$. As SO(2) acts by isometries,

$$d^2(x, y) = \inf_\theta \sum_i |x_i - e^{i\theta} y_i|^2$$
$$= \inf_\theta \left[|x|^2 + |y|^2 - 2x \cdot y \cos\theta + 2\omega(x, y) \sin\theta\right].$$

Taking the derivative with respect to $\theta$, we see that the minimum occurs at $\theta = \theta_0$ where $x \cdot y \sin\theta_0 + \omega(x,y) \cos\theta_0 = 0$. This implies that

$$d^2(x,y) = |x|^2 + |y|^2 - 2\sqrt{(x \cdot y)^2 + \omega(x,y)^2} = |x|^2 + |y|^2 - 2|\langle x, y \rangle|.$$

The $\varepsilon^2$ term in the expansion of $d^2(x, x + \varepsilon v)$ is the reduced kinetic energy $K_{\text{red}}(x,v)$ corresponding to the decomposition $K_{\text{red}} = K - K_{\text{rot}}$. We find

$$K_{\text{red}}(x,v) = |v|^2 - \frac{\omega(x,v)^2}{|x|^2}.$$

(This is the expression of the pullback of the natural induced metric on the quotient $\mathcal{X}/SO(2)$ by the projection of $\mathcal{X}$ onto the quotient.) This is consistent since $K_{\text{red}}(x,v) = 0$ when $v$ is tangent to the orbit of $x$, i.e. when $v$ is proportional to $ix$.

(iii) *The length $\ell_0$ in spherical coordinates.* It will be convenient to use spherical coordinates defined in the shape space $\mathbb{R}^3$ by

$$u_1 = r^2 \cos\varphi \cos\theta, \; u_2 = r^2 \cos\varphi \sin\theta, \; u_3 = r^2 \sin\varphi.$$

We have

$$r^2 = \sqrt{u_1^2 + u_2^2 + u_3^2} = I,$$

which justifies the choice of the notation $r$. A tedious but simple calculation, or an appeal to the $U(2)$-invariance of $K_{\text{red}}$, now proves

LEMMA 5 (see [7]). *In spherical coordinates the quotient metric corresponding to the reduced kinetic energy $K_{\text{red}}$ occurring in the reduced action is given by*

$$ds^2 = dr^2 + \frac{r^2}{4}(\cos^2\varphi \, d\theta^2 + d\varphi^2).$$

*In particular the shape sphere $I = r^2 = 1$ is isometric to the standard sphere of radius $1/2$, and the shape space $\mathbb{R}^3$ is the cone over this sphere, the sphere itself consisting of all those points at distance 1 from triple collision.*

Using Lemma 4 we can express the equipotential curve through the Euler configurations on the shape sphere by the implicit equation:

$$\frac{1}{\sqrt{1 + \cos\varphi \cos\theta}} + \frac{1}{\sqrt{1 + \cos\varphi \cos(\theta + \frac{2\pi}{3})}} + \frac{1}{\sqrt{1 + \cos\varphi \cos(\theta + \frac{4\pi}{3})}} = \tilde{U}_E = \frac{5}{\sqrt{2}}.$$

This curve is a double covering of the equator $\varphi = 0$ and as such may be parametrized by a function $\varphi = \varphi(\theta)$, provided $\theta$ is allowed to vary in an interval of length $4\pi$. We are interested in $\ell_0$ which is one twelfth of its length. It follows from Lemma 5 that

$$\ell_0 = \frac{1}{12} \int_0^{2\pi} \sqrt{\cos^2 \varphi(\theta) + \varphi'^2(\theta)} \ d\theta = \frac{1}{2} \int_0^{\frac{\pi}{3}} \sqrt{\cos^2 \varphi(\theta) + \varphi'^2(\theta)} \ d\theta.$$

To finish the proof of Lemma 3, and hence of the existence of the collision free minimizer, we use the following numerical estimate of the Euler equipotential length $\ell_0$ obtained by Carles Simó and later confirmed by Jacques Laskar:

$$\frac{\pi}{5.082553924511} \leq \ell_0 \leq \frac{\pi}{5.082553924509}.$$

These estimates were obtained by using a Newton method for computing $\varphi(\theta)$ and $\varphi'(\theta)$, and then the trapezoid method for computing the integral.

*Remark.* We explain the meaning of the spherical coordinates in terms of triangles. The parallels or "latitudes" $\varphi = $ constant in the shape sphere correspond to triangles with the same orientation and the same ellipse of inertia up to rotation. Indeed, this set of triangles is characterized by a common area (see [1]). But the area is proportional to $\mathrm{Im}\bar{z}_1 z_2$, that is to $u_3 = \sin\varphi$, the height function on the sphere. The meridians or "longitudes" $\theta = $ constant in the shape sphere are defined by a linear relation between the squares of the mutual distances. These properties of the coordinates $(\theta, \varphi)$ are a consequence of the invariance of the metric under the orthogonal group $\mathrm{O}(\mathcal{D})$ of the disposition space $\mathcal{D}$ (see [1] for the definition). The equilateral triangles $L_\pm$ (north and south poles) are fixed points of the action of $\mathrm{SO}(\mathcal{D})$ and the parallels above are the circles with center $L_\pm$ which are also the orbits of the action of $\mathrm{SO}(\mathcal{D})$. The geodesics through $L_\pm$ orthogonal to these circles are the meridians. They are transitively transformed into each other by $\mathrm{SO}(\mathcal{D})$ and each one is fixed by an involution in $\mathrm{O}(\mathcal{D})$.

## 6. Symmetries: Proof of the existence of the "eight"

Recall our action of the Klein group $\mathbb{Z}/2\mathbb{Z} \times \mathbb{Z}/2\mathbb{Z}$ on $\mathbb{R}/\overline{T}\mathbb{Z}$ and $\mathbb{R}^2$. If $\sigma$ and $\tau$ are generators,

$$\sigma \cdot t = t + \frac{\overline{T}}{2}, \ \tau \cdot t = -t + \frac{\overline{T}}{2}, \quad \sigma \cdot (x,y) = (-x, y), \ \tau \cdot (x, y) = (x, -y).$$

LEMMA 6. *After being symmetrized according to the pattern of the Euler equipotential curve on the shape sphere, a minimizing path $x$ gives a $\overline{T} = 12T$ periodic loop (still called $x$), with zero angular momentum, which, up to a time translation and a space rotation, is of the form*

$$x(t) = \bigl(q(t + 2\overline{T}/3), q(t + \overline{T}/3), q(t)\bigr),$$

*described in the theorem.*

The proof of Lemma 6 proceeds in three steps.

*Step* 1. Observe that the minimizing arc is orthogonal to the two manifolds $E_3$ and $M_1$ constraining its endpoints. This follows from the boundary term arising in the first variation formula for the action.

*Step* 2. Observe that upon reflecting this arc about one of the three meridians, or about the equator, we will obtain another minimizing solution arc, one with permuted endpoint conditions; *e.g.* with $E_j$ and $M_k$ in place of $E_3$ and $M_1$. Using these reflections we build the entire closed solution curve in the reduced configuration space. It consists of 12 subarcs all congruent to our original minimizer. Orthogonality guarantees that they fit together smoothly, thus forming a single solution.

More precisely, because reflection about the meridian $M_1$ is a symmetry of the reduced action (and hence of the equations), and because the minimizing arc is orthogonal to the meridian at its endpoint, when we continue the solution represented by the arc through the meridian $M_1$, the result is the same as if we had reflected it about the meridian, and then reversed the direction of time. In symbols: $x(\overline{T}/12 + t) = s_1(x(\overline{T}/12 - t))$ where $s_1$ is reflection about the meridian $M_1$, and where $\overline{T} = 12T$ will be the period of the full orbit. ($T$ is the time it takes to hit the meridian.)

The reflection $s_1$ can be realized in the inertial plane as follows: at time $T = \overline{T}/12$ the triangle is an isosceles triangle of type $M_1$ and hence has a reflectional symmetry $\tau$. Choose coordinates in $\mathbb{R}^2$ so that $\tau(x,y) = (x, -y)$, i.e. so that the perpendicular bisector of the edge joining 2 to 3 is the $x$-axis. Let $S_1(x_1, x_2, x_3) = (x_1, x_3, x_2)$ be the operation of interchanging masses 2 and 3. Then $s_1 = S_1 \circ \tau = \tau \circ S_1$.

We now have a solution from $E_3$ to $E_2$ in time $2T = \overline{T}/6$. By a similar argument, to continue this arc of solution through $E_2$ we must perform a half-twist $H_2$ through $E_2$ and reverse time: $x(2T - t) = H_2(x(2T + t))$. This half-twist is a symmetry of the action being the composition of reflection about the equator with reflection about the meridian $M_2$. It is realized in the inertial plane by interchanging masses 1 and 3 and then performing the inertial half twist $\sigma \circ \tau(x, y) = (-x, -y)$ about the origin.

We obtain in this way an arc of solution from $E_3$ to $E_1$ in time $4T = \overline{T}/3$.

Continuing around the equator in this manner with the appropriate reflections or half twists we construct a smooth curve in the reduced configuration space which consists of 12 congruent arcs, alternating in pairs above and below the equator, so as to have the same symmetry as the equipotential curve. It is a solution curve, and is $\overline{T}$-periodic mod rotations.

*Step* 3. We have constructed the projection of our solution curve to the reduced configuration space. We now reconstruct the full solution curve, show that it is periodic (i.e. in inertial space), and show that it satisfies all the properties described in the theorem. This is done by invoking the area rule for reconstructing the original dynamics from the reduced dynamics, and by using the symmetries of the curve.

Figure 3 shows segments of the reduced orbit on the shape sphere and anticipates the reconstructed orbit in the inertial plane.

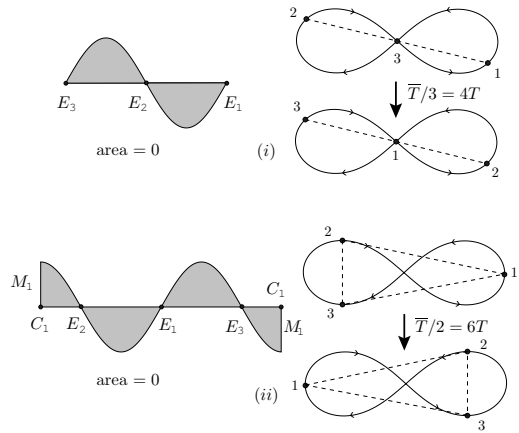

Figure 3

*The area rule* tells us how to recover the motion of the masses in the inertial plane, given the curve representing this motion in the shape sphere. Suppose the shape curve is closed. Then the initial and final triangles in the plane are similar. The angle of rotation which relates these two triangles up to scale equals twice the spherical area enclosed by the shape curve. (The area of the sphere of radius $1/2$ is $\pi$. The factor of 2 in the area formula insures that the answer is well-defined modulo $2\pi$.) For a proof of the area rule see, for example, [9] or references therein. If instead the shape curve is not closed, but begins and ends on the equator of collinear shapes, then we close it up by travelling "backwards" along the equator from the endpoint until we reach the beginning point. Now compute twice the signed area enclosed by this closed curve. This equals the angle between the two lines in the inertial plane which contain the initial and final configuration for any zero angular momentum curve realizing the given shape curve*. Finally, if the curve begins or ends on one of our three isosceles meridians, we compute the angle between the beginning and ending symmetry axis of the isosceles triangle by following the appropriate meridian up or down to the equator, travelling along the equator to close up the curve, and then computing the area within the resulting closed curve.

---

*There are two ways to close up the curve into a loop, depending on which way we travel the equator. The two angles so computed differ by $\pi$, this being twice the area of a hemisphere of a sphere of radius one-half. This is no problem since the angle between two unoriented lines is only defined modulo $\pi$.

The fact that the signed areas depicted on this figure are equal to zero implies that as we travel our solution curve

(i) if we start at an Euler configuration and follow the orbit for a time $\overline{T}/3 = 4T$, passing through an intermediate Euler configuration at time $2T$ we arrive at an Euler configuration with the three masses sitting on the same line as that of the initial Euler configuration. That is to say, there is no rotation of the Euler line, contrary to what happened at the intermediate time $2T$.

(ii) after time $\overline{T}/2$, an isosceles configuration returns to itself, up to reflection (equivalently up to interchange of the symmetric vertices). There is no rotation of the symmetry axis of the triangle.

Choose the origin of time $t = 0$ to correspond to being in the Euler configuration $E_3$. Set $q(t) = x_3(t)$ where $x(t) = (x_1(t), x_2(t), x_3(t))$ is our solution. The first property implies that

$$\begin{cases} q(t) = x_3(t) \text{ for } 0 \leq t \leq \overline{T}/3, \\ q(t) = x_2(t - \overline{T}/3) \text{ for } \overline{T}/3 \leq t \leq 2\overline{T}/3, \\ q(t) = x_1(t - 2\overline{T}/3) \text{ for } 2\overline{T}/3 \leq t \leq \overline{T}. \end{cases}$$

Thus, after time $\overline{T}/3$ (respectively, $2\overline{T}/3$), the bodies $2, 3, 1$ have been replaced by the bodies $3, 1, 2$ (respectively, $1, 2, 3$) with the same velocities. The three bodies move along the same closed curve $q(t)$ of period $\overline{T}$ with a phase shift relative to each other of $\overline{T}/3$. Using (ii), the Klein symmetry follows easily.

Figure 4 shows the projection (still called $x$) of the orbit in the reduced configuration space.

*Step* 4. It remains to be proven that the equivariant curve $q$ we constructed not only has the required symmetry but also has the shape of a figure eight without extra small loops or other unpleasant features. We will use the following basic lemma:

LEMMA 7. *The angular momentum $q(t) \wedge \dot{q}(t)$ is nowhere vanishing for $t < 0 \leq \overline{T}/4$; i.e. in any one-quarter of the curve the angular momentum of the mass tracing out that quarter is nonzero.*

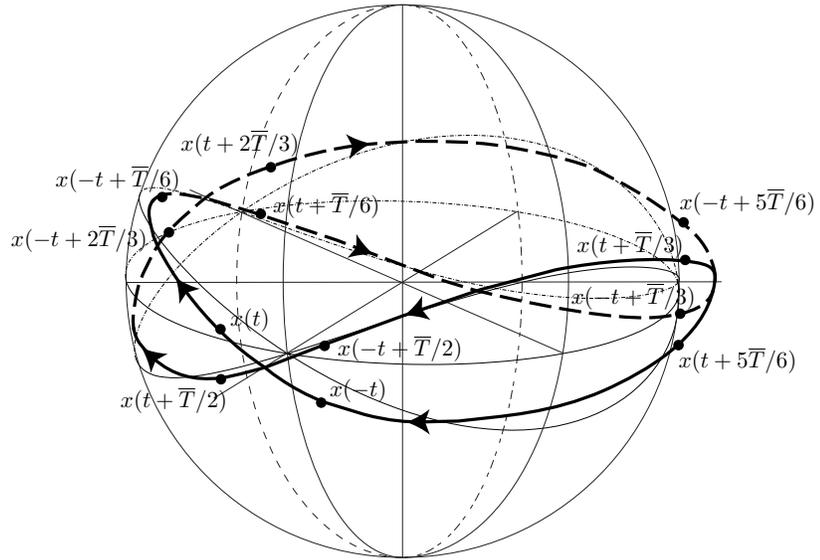

Figure 4

*Proof of Lemma* 7. We will need the fact that Newton's equations are satisfied, and two consequences of the minimality of the curve in the shape sphere.

We begin by computing the time derivative of $q \wedge \dot{q}$. It is

$$\frac{d}{dt}(q \wedge \dot{q}) = q \wedge \frac{d^2 q}{dt^2}.$$

Now use the fact that $q(t) = x_3(t)$ where $(x_1(t), x_2(t), x_3(t))$ satisfy Newton's equations. Also use the fact that the center of mass is zero at all times: $x_1 + x_2 + x_3 = 0$. These yield

$$\frac{d}{dt}(q \wedge \dot{q}) = \left(\frac{1}{r_{13}^3} - \frac{1}{r_{23}^3}\right)(x_3 \wedge x_1).$$

There are only two ways this quantity can be zero: either

(A) $x_1$ and $x_3$ are linearly dependent, or

(B) $r_{13} = r_{23}$.

To eliminate these possibilities, we use the reflection principle in shape space. First, if $x_1$ and $x_3$ are linearly independent, then the entire configuration is collinear. Thus it crosses the equator in shape space. This can happen for a minimizing arc only at an Euler point. For if it did at another point, the minimizing arc between Euler and Isosceles ($0 < t < \overline{T}/12$) would be divided into two (or more) subarcs, one of which lies below the equator and the other above. Reflect one of these arcs across the equator, leaving the other(s) fixed.

The resulting arc has the same action, has no collisions, and still connects Euler to Isosceles, and hence is also a collision-free minimizer. But it is no longer analytic, contradicting the fact that collision-free minimizers correspond to solutions.

Next suppose that $r_{13}(t) = r_{23}(t)$ at some time $t$, with $0 < t < \overline{T}/12$. This asserts that the curve has recrossed the isosceles meridian passing through our initial Euler point. The same reflection principle applies, with the Euler meridian playing the role of the equator.

We can now suppose that, as in Figure 4, our minimizing arc lies in the left upper quarter of the shape sphere for $0 < t < \overline{T}/6$. The arguments thus far show that the angular momentum $q(t) \wedge \dot{q}(t)$ decreases from the value $0$ for $0 < t < \overline{T}/6$ and increases for $\overline{T}/6 < t < \overline{T}/4$ (a bivector in the plane is called positive if it is a positive multiple of the standard area form $e_1 \wedge e_2$). It remains to notice that it takes a negative value at $t = \overline{T}/4$ to conclude that it stays strictly negative for $0 < t < \overline{T}/4$. It follows that it stays strictly negative for $0 < t < \overline{T}/2$ and strictly positive for $\overline{T}/2 < t < \overline{T}$.

COROLLARY. *Each lobe of the curve is starshaped: the only time $x_i \wedge \dot{x}_i$ becomes zero, for $i = 1, 2, 3$, is when $x_i$ passes through the origin.*

*Proof of the corollary.* In polar coordinates $(r, \theta)$ the angular momentum has the well-known expression $x \wedge \dot{x} = (r^2 \dot{\theta}) e_1 \wedge e_2$. From this it follows that the polar angle $\theta(t)$ of the curve $q(t)$ decreases monotonically over the interval $(0, \overline{T}/2)$ from its maximum value of $\theta(0)$ when $r(0) = 0$ to its minimal value of $\theta(\overline{T}/2) = -\theta(0)$.

## Appendix 1. Estimates

We get estimates on $U$ and $K$ along our solution. These imply estimates on $I$ and $J$ by using the fact that the solution has zero mean, which follows from its symmetries.

Let

$$I_0 = \left(\frac{5}{\sqrt{2}\ell_0^2}\right)^{\frac{2}{3}} T^{\frac{4}{3}}, \quad U_0 = K_0 = \frac{\ell_0^2 I_0}{T^2},$$

be defined as above and

$$H_0 = \frac{1}{2}K_0 - U_0 = -\frac{1}{2}U_0, \quad a = (\frac{1}{2}K_0 + U_0)T = \frac{3}{2}U_0 T = -3H_0 T,$$
$$J_0 = \sqrt{I_0 K_0}.$$

Let $\langle f \rangle = \frac{1}{T} \int_0^T f(t)\, dt$ be the mean value of a function $[0, T] \to R$.

LEMMA 8. *A minimizing path satisfies the following estimates*:

$$\langle U \rangle = \langle K \rangle < U_0 = K_0, \quad H > H_0, \quad \langle I \rangle < \frac{36\ell_0^2}{\pi^2} I_0, \quad \langle |J| \rangle < \frac{6\ell_0}{\pi} J_0.$$

*Proof of Lemma 8.* Because the conservation of energy is true almost everywhere on a minimizing path*, such a path $x$ has action

$$\mathcal{A}(x) = HT + 2\langle U \rangle T < a = H_0 T + 2U_0 T.$$

But we deduce from the Lagrange-Jacobi identity $\dot{J} = \frac{1}{2}\ddot{I} = 2H + U = K - U$ that

$$\langle K \rangle - \langle U \rangle = J(T) - J(0).$$

Because of the symmetry of the minimizing path, $J(T) = J(0) = 0$ (note that the eventual presence of a double collision at the end of the path would not change this fact because then $J(t)$ would behave as $(T-t)^{\frac{1}{3}}$ which goes to 0 as $t$ goes to $T$). This implies

$$\langle K \rangle = \langle U \rangle = -2H,$$

so that

$$\mathcal{A}(x) = -3HT < a = -3H_0T \quad, \text{ hence } \quad H > H_0.$$

The inequalities concerning $\langle U \rangle$ and $\langle K \rangle$ follow immediately.

To bound $\langle I \rangle$, note that by construction our $x$ has zero average (in $\mathcal{X}$) over its full period $12T$. By the Poincaré lemma this implies $\langle K \rangle > \frac{4\pi^2}{(12T)^2} \langle I \rangle$ and consequently we get the estimate on $\langle I \rangle$.

Finally, Sundman's inequality $IK - J^2 \geq 0$ yields the bound on $\langle |J| \rangle$.

*Simó's numerical computation of actions.* The following numerical estimates of various important actions were obtained by C. Simó. The period here is taken to be $T = 2\pi/12$.

$$A_2 = (\frac{3}{2})^{2/3}\frac{\pi}{2} = 2.0583255\ldots,$$

$$a = \mathcal{A}(\text{test}) = (225\pi l_0^2/32)^{1/3} = 2.0359863\ldots,$$

$$A_{\min} = \mathcal{A}(\text{solution}) = 2.0309938\ldots.$$

It is interesting to notice that the action of any element of $\Lambda$ undergoing a triple collision is much higher. Indeed, Sundman's inequality implies it is higher than the action $A_3$ of an equilateral homothetic solution from collision to zero velocity. This last action is given by the same formula as $A_2$ with the value $\tilde{U}_3 = 3$ of $\tilde{U}$ at the equilateral configuration replacing $\tilde{U}_2 = 1/\sqrt{2}$, so that

$$A_3 = \left(3\sqrt{2}\right)^{\frac{2}{3}} \times A_2 = 5.39433\ldots.$$

---

*This classical result can be proved by computing the variation of action due to a change of parametrization.

## Appendix 2. How this orbit was discovered

One of us (R.M.) had been searching for several years for periodic orbits in the three-body problem using the method of minimization of the action over well-chosen homotopy classes of loops in the configuration space. (See [7] where the collision problem is avoided by a strong force hypothesis, and where potentially interesting homotopy classes are described.) Approximately at the same time, but independently, both authors realized that equality among the masses could make the variational approach more tractable by allowing us to impose additional symmetries on competing loops. This led to the preprint [3] by A.C. and A. Venturelli, submitted to *Celestial Mechanics*, in which new periodic orbits were found for the spatial four body problem with equal masses, and to the preprint [8], submitted by R.M. to Nonlinearity. A.C. was asked to act as a referee for preprint [8], titled "Figure eights with three bodies", which described two different types of periodic orbits of the three-body problem according to whether or not the masses were all equal or only two were equal. Only the orbit for all equal masses survived careful scrutiny by A.C. and A. Venturelli. The other orbit – which, curiously, was the one which had given the paper its title – was supposed to be a figure eight not in the plane, but in the shape sphere. However the proof of the absence of collisions for this orbit was found to be in error. In trying to understand the case of equal masses, A.C. discovered, at first experimentally and then mathematically, that the three equal masses must travel along a fixed eight-shaped curve in the plane. This discovery gave new life to the title of the preprint. The numerical experiment grew out of an example called "figure eight attractor" in the nice program "Gravitation" by Jeff Rommereide. The success of the experiment came from the constraints placed on the velocities $v$ at any Euler point of the orbit. These constraints, depicted on figure 1 for $E_3$ as $v_1 = v_2 = -v_3/2$, are due to the fact that the angular momentum is zero and that each Euler configuration is an extremum of the moment of inertia $I$ along the orbit: $dI(v) = 0$. By making more stringent symmetry assumptions on the orbit than R.M had made, A.C. was then able to give a direct and simple proof, partly in the spirit of [2], of the absence of triple collisions, and to obtain estimates for $I$ and $U$ along the orbit. Finally, R.M. noticed the trick of forgetting one mass, which made the calculations for triple collisions extend to double collisions and bypassed completely any local variational analysis. Precise numerical computations by Carles Simó, using the special form of velocities at an Euler configuration, gave accurate pictures of the eight and showed in particular its "stability"*.

*Acknowledgements.* The authors thank warmly Andrea Venturelli and Carles Simó for many enlightening discussions, and again Carles Simó for his crucial numerical help. Warm thanks also from A.C. to Jean Petitot who, some years ago, gave him the program "Gravitation". R.M. thanks Neil Balmforth for help visualizing the eight, Julian Barbour and Wu-Yi Hsiang for inspirational conversations. Thanks to Jacques Laskar who, as an editor of Nonlinearity, authorized the referee of [8] to contact the author. And finally, thanks again to Jacques Laskar for producing the final form of the figures.

Astronomie et Systèmes Dynamiques, IMCCE, UMR 8028 du CNRS, Paris, France
*E-mail address*: chencine@bdl.fr

Université Paris VII-Denis Diderot, Paris, France

UCSC, Santa Cruz, CA
*E-mail address*: rmont@math.ucsc.edu

---

*This is particularly interesting because very few stable periodic orbits in the inertial frame are known for the three-body problem with equal masses. One example is Schubart's collinear orbit: Nulerische Aufsuchung periodischer Lösungen im Dreikörperproblem, Astronomische Nachriften vol. 283, pp. 17–22, 1956 (thanks to C. Simó for this reference). One fourth of this orbit travels collinearly from $E_3$ to $C_1$. The "eight" has strong similarities with the critical point which would result from Schubart's orbit by permuting the colliding bodies at each collision so that the result travels completely around the equator in the shape sphere. Another example is that of orbits periodic in a rotating frame in case of resonance. To this category belong the tight binary solutions with a third mass far away and the family of "interplay" solutions connecting them to Schubart's orbit (see M. Hénon, A family of periodic solutions of the planar three-body problem, and their stability, Celestial Mechanics 13, pp. 267–285, 1976 and the references there to papers by R. Broucke and J. D. Hadjidemetriou).

*Note added in proof.* After this paper was accepted for publication, Phil Holmes brought the work of C. Moore [10] to our attention. Moore had already found the figure eight solution numerically, using gradient flow for the action functional.

Also, in [11], K.-C. Chen proves a better estimate than $A_2$ for loops in $\Lambda$ with collisions. This allows him to replace the equipotential test path by pieces of great circles and to conclude without using a computer.